\documentclass[10pt]{amsart} 

\usepackage{amsfonts,amsmath,amsthm, graphics}

\title[R.F. Shamoyan, E.B. Tomashevskaya]{On Bergman type projection in some new analytic function spaces in bounded strongly pseudoconvex domains}

\author[]{R.F. Shamoyan, E.B. Tomashevskaya}

\begin{document}
\begin{abstract} We prove the boundedness of Bergman type projections in two different analytic function spaces in bounded strongly pseudoconvex domains with the smooth boundary. Our results were previously well-known in the case of the unit disk. \\
Keywords: pseudoconvex domains, unit disk, analytic function, Bergman projection.
\end{abstract}

\maketitle 

\section{INTRODUCTION}

We refer the interested reader to the classical and detailed definition of bounded strongly pseudoconvex domains   in $\mathbb C^n$ via defining function, Hessian of a defining function and a complex tangent space to [4].

Let  $G\subset {\mathbb{C}}^n$ be a domain, that is, an open connected subset. One says that $G$ is pseudoconvex (or Hartogs pseudoconvex) if there exists a continuous plurisubharmonic function $\varphi$ on $G$ such that the set  $\{ z \in G \mid \varphi(z) < x \}$ is a relatively compact subset of $G$ for all real numbers $x.$ In other words, a domain is pseudoconvex if $G$ has a continuous plurisubharmonic exhaustion function. Every (geometrically) convex set is pseudoconvex. However, there are pseudoconvex domains which are not geometrically convex.

When $G$ has a $C^2$ (twice continuously differentiable) boundary, this notion is the same as Levi pseudoconvexity, which is easier to work with. More specifically, with a $C^2$ boundary, it can be shown that  $G$ has a defining function, i.e., that there exists  $\rho: \mathbb{C}^n \to  \mathbb{R}$ which is $C^2$ so that $G=\{\rho <0 \},$ and $\partial G =\{\rho =0\}.$  Now,  $G$ is pseudoconvex iff for every  $p \in \partial G$ and $w$ in the complex tangent space at $p,$ that is,  

\[\nabla \rho(p) w = \sum_{i=1}^n \frac{\partial \rho (p)}{ \partial z_j }w_j =0,\] we have
\[\sum_{i,j=1}^n \frac{\partial^2 \rho(p)}{\partial z_i \partial \bar{z_j} } w_i \bar{w_j} \geq 0.\]

The definition above is analogous to definitions of convexity in Real Analysis.

We provide an extension of a classical result to pseudoconvex domains. 
Throughout this paper $H(D)$ or $H(\Lambda)$ denotes the space of all holomorphic functions on an open set $D\subset C^n$ or  $\Lambda\subset C^n$. We denote by dv the normalized Lebegues measure on $D$.

We follow notation from [5]. Let $D$ be a bounded strongly pseudoconvex domain in $C^n$ with smooth boundary, let $d(z)=dist(z, \partial  D)$. Then there is a neighbourhood $U$ of $\bar{D}$ and $\rho \in C^\infty(U)$ such that $D=\{z\in U:\rho(z)>0\}$ , $ |d\rho(z)| \geq c> 0$ for $z\in \partial D, 0<\rho (z)<1$ for $z\in D$ and $-\rho$ is strictly plurisubharmonic in a neighborhood  $U_0$ of $\partial D$. Note that $d(z) \asymp \rho(z), z\in D$.

Then there is an $r_0$ such that the domains $D_r= \{z\in D : \rho(z)>r\}$  are also smoothly bounded strictly pseudoconvex domains for all $0<r\geq r_0$. Let $d\sigma_r$ is the normalized surface measure on $\partial D_r.$
The following analytic mixed norm spaces were investigate in [5]. 

For $0<p<\infty, 0<q\leq \infty, \delta>0, k=0,1, r$ we set
$$
||f||_{p,q,\delta, k}=\sum\limits_{|\alpha|
\leq k}\left(\int\limits_{0}^{r_0}\left(r^\delta\int\limits_{\partial D_r}
|D^\alpha f|^pd\sigma_r\right)^{q/p}\frac{dr}{r}\right)
^{1/q},  0<q<\infty;
$$

$$
||f||_{p,\infty,\delta, k}=\sup\limits_{0<r<r_0}  \sum\limits_{|\alpha|
\leq k}\left(r^\delta\int\limits_{\partial D_r}
|D^\alpha f|^pd\sigma_r\right)^{1/p}.
$$

The corresponding $A^{p,q}_{\delta,k}=\{f\in H(D): ||f||_{p,q,\delta, k}<\infty\}$ are complete quasinormed spaces for $p,q\geq 1.$  They are Banach spaces. We mostly deal with the case $k=0$ then we write simply $A^{p,q}_{\delta}$ and  $||f||_{p,q,\delta}$.
  
 We also consider these spaces for $p=\infty$ and $k=0$, the corresponding space is denote by  $A^{\infty,p}_{\delta}=A^{\infty,p}_{\delta}(D)$ and consists of all $f\in H(D)$ such that 

$$
||f||_{\infty,p ,\delta}=\left(\int\limits_0^{r_0}\left(\sup\limits_{\partial D_r}|f|\right)^p\cdot r^{\delta p-1}dr\right)^{1/p}<\infty.
$$

Also for $\delta >-1$ the space $A^{\infty}_{\delta}=A^{\infty}_{\delta}(D)$ consists of all $f\in H(D)$ such that 

$$
||f||_{A_\delta ^\infty}=\left(\sup\limits_{z\in D}\right)|f(z)|\rho(z)^\delta<\infty.
$$

and the weighted Bergman space $A^p_\delta$.

$A^p_\delta=A^{p,p}_{\delta+1}(D)$ consists of all  $f\in H(D)$ such that 
$$
||f||_{A_\delta ^p}=\left(\int\limits_{D}|f(z)|^p\rho(z)^\delta dv(z)\right)^{1/p}<\infty,
$$
where $dv$ is a normalized Lebegues measure in $D$.

Let further $dv_\alpha=(\delta^\alpha)dv(z)$ where $\delta(z)=\rho(z).$ Our proofs are heavily based on the estimates from [3] where more general situation was considered. Since $|f(x)|^p$ is subharmonic (even plurisubharmonic) for a holomorphic $f$ we have $A^p_s(D)\subset A^\infty_t(D)$, for $0<p<\infty, sp>n$ and $t=s$
also  $A^p_s(D)\subset A^1_s(D)$ for $0<p\leq 1$ and $A^p_s(D)\subset A^1_t(D)$ for $p>1$ and $t$
sufficiently large.

Therefore we have an integral representation
$$
f(z)=\int\limits_D f(\xi) K_{\tilde{t}}(z,\xi)\rho^t(\xi)dv(\xi),  \eqno(1)
$$
$f\in A^1_t(D), z\in D, \tilde{t}=t+n+1$,
where $ K_{\tilde{t}}(z,\xi)$ is a kernel of type $t$ that is a smooth function on $D\times D$ such that 
$$
| K_{\tilde{t}}(z,\xi)| \leq C_1|\tilde{\Phi}(z,\xi)|^{-(n+1+t)},
$$
where $\tilde{\Phi}(z,\xi)$ is so called Henkin-Ramirez function on $D$ (see for detailed explanations of these issues related with the Bergman kernel and integral representations in pseudoconvex domains  for example [3], [5] and references there).

Note that (1)  holds for function in any space $X$ that embeds into $A^1_t$. 
We review some known facts on $\tilde{\Phi}$.
 This function is $C^\infty$ in $U\times U$ where $U$ is a neighborhood 
 of $\bar{D}$ it is holomorphic in $z$ 
 and $\tilde{\Phi}(\xi,\xi)=\rho(\xi)$ for $\xi \in U$. Moreover on $\bar{D}\times \bar{D}$ it vanished only on the diagonal $(\xi,\xi), \xi \in \partial D$. 

Locally it is up to a non vanishing smooth multiplicative factor equal to the Levi polynomial of $\rho$.

We denote various positive constants in this paper by $C, C_1, C_2$  ect.

{\bf Lemma 1} (see [1]-[3], [5]). 

Assume $K_{\tilde{t}}(z,\xi)$ is a kernel of type  $t; t>-1, \tilde{t}=t+n+1.$ 

a) For $0<r<r_0$ we have 
$$
\int\limits_{\partial D_r}|K_{\tilde{t}}(z,\xi)|d \sigma_r(z)\leq C_2(\rho(\xi)+r)^{-t-1}; \xi\in D.
$$

b) Assume $\sigma>0$ satisfies $\sigma-t-1<0.$ Then we have

$$
\int\limits_{D}|K_{\tilde{t}}(z,\xi)|\rho^{\sigma-1}(z)dv(z)\leq C_3\rho^{\sigma-t-1}(\xi); \xi\in D.
$$

Note that the same estimates are valid if $K$ is replaced by $\tilde{K}_{t}(z,\xi)=K_t(\xi,z), t>0.$

Some facts on $A^{p,q}_{\alpha,k}$ space.

{\bf Proposition 1} (see [1]-[3], [5]).

If $0<p_0<p_1<\infty, 0<q\leq \infty$, $\delta_0, \delta_0'>0$ and $\frac{n+\delta_0'}{p_1}=\frac{n+\delta_0}{p_0}$ then

$A^{p_0,q}_{\delta_0,k}(D)\subset A^{p_1,q}_{\delta_0',k}(D)$.

{\bf Corollary 1}  (see [1]-[3], [5]).

If $0<p\leq 1, \alpha>-1$ then 
$$A_{\alpha}^{p}(D)\subset A_{\beta}^{1}(D), \beta=\frac{n+1+\alpha}{p}-(n+1),$$
that is 
$$
\int\limits_D|F(z)|\rho (z)^{\frac{n+1+\alpha}{p}-(n+1)}dv(z)\leq C_4||F||_{A^p_\alpha}, F\in H(D).
$$

{\bf Proposition 2} (see [1]-[3], [5]). 

If $0<p_0<\infty, 0<q_0< q_1\leq \infty,$ $\delta_0>0$ and $k=0,1,2,$ then 
$A^{p,q_0}_{\delta,k}(D)\subset A^{p,q_1}_{\delta,k}(D)$.

{\bf Lemma 2.} 

Let $\delta>0, t>t_0$, $t_0$ is large enough then, $r_0$ is fixed 
$$
r^\delta\cdot\left(\int\limits^{r_0}_{0}\frac{R^{t-\delta}dR}{(r+R)^{t+1}}\right)\leq C_5.
$$
 
{\bf Lemma 3}  (see [4]). 
 
Let  $0<p<1,s>-1,$ $ r>0, $ $t=p(s+n+1)-(n+1).$ Then we have (even if we replace $|\tilde{\Phi}|^r$ by $K$, where $K$ is a kernel of $r>$ type.
$$
\left| \int\limits_D |f(\xi)|\cdot |\tilde{\Phi}(z,\xi)|^r\cdot d(\xi)^sdv(\xi)\right|^	p \leq  C_6\int\limits_D |f(\xi)|^p\cdot |\tilde{\Phi}(z,\xi)|^{rp}\cdot  d(\xi)^tdv(\xi).
$$

{\bf Remark 1.}  This Lemma is valid if we replace $\alpha$ by $\rho$ since $d(\xi)\asymp \rho(\xi)$. Below we denote this function also by $\delta (\xi).$ 

\section{MAIN RESULTS}

In this section we formulate our main results.

{\bf Theorem 1.} 
 
Let $0<p, q\leq 1,$ $\alpha>\tilde{\alpha_0}),$ $\tilde{\alpha_0}$ is large enough then
$$
\left(T_\alpha f\right)(w)=\int\limits_\Lambda K_{\alpha_0}(z,w)f(z) dv_\alpha(z), \alpha_0=\alpha+n+1,
$$
maps $A_\beta^{p,q}$ into  $A_\beta^{p,q}$ for all $\beta, \beta>0$, $dv_\alpha=\delta^\alpha dv.$

{\bf Remark 2.}  This theorem in the unit disk is a well known classical result about Bergman projection in the unit disk.

{\bf Proof of theorem 1.}  By lemma 3 above 
$$
\left|\left(T_\alpha f\right)(w)\right|^p\leq C_7\int\limits_\Lambda |K_{\alpha_0}(z,w)|^p(\delta(z))^{(n+1)(p-1)+(\alpha)p}\cdot|f(z)|^p\cdot dv(z). 
$$

Then by lemma 1 we have
$$
\int\limits_{\partial\Lambda_{\tilde{r}}}\left|\left(T_\alpha f\right)(w)\right|^p d\sigma_{\tilde{r}}\leq C_8\int\limits_0^{r_0}\int\limits_{d\Lambda_r} \left(\delta(z)^{\alpha p+(n+1)(p-1)}\right)\cdot|f(z)|^p\cdot\frac{1}{(r+\delta(z))^{(\alpha+n+1)p-(n)}}\cdot w(z)d\sigma _\varepsilon dr, 
$$
$w(z)\in C^\infty\Lambda)$.

Let $q\leq p.$ Then we have that using lemmas above and the fact that $\int\limits_{\partial\Lambda_r}|f|^p$ is monotone. 

$$
\int\limits_{0}^{\tilde{r_0}}\tilde{r}^{q/p\beta}\left(\int\limits_{\partial\Lambda_{\tilde{r}}}\left|\left(T_\alpha f\right)(w)\right|^p d\sigma_{\tilde{r}}\right)^{q/p}\times \left(\frac{d\tilde{r}}{\tilde{r}}\right)\leq  
$$

$$
\leq C_9\int\limits_0^{\tilde{r}_0}\left(\tilde{r}^{\beta q/p-1}\right) \times
$$

$$\times\left(\int\limits_{0}^{r_0}\int\limits_{d\Lambda_r} \left(\delta(z)^{\alpha p+(n+1)(p-1)}\right)\cdot|f(z)|^p\cdot\frac{1}{(r+\delta(z))^{(\alpha+n+1)p-(n)}}\cdot d \sigma _\varepsilon dr\right)^{q/p} \cdot d\tilde{r}\leq
$$

$$
\leq C_{10}\int\limits_0^{\tilde{r}_0}\left(\tilde{r}^{\alpha q+(n+1)q-nq/p-1}\right) \times
$$

$$\times\left(\int\limits_{d\Lambda_r} |f(z)|^p d\sigma\right)^{q/p}\cdot\left(\int\limits_{0}^{\tilde{r}_0}\frac{\tilde{r}^{\beta q/p-1}\cdot d\tilde{r}}{(r+\tilde{r})^{(\alpha+n+1)q-nq/p}} \right)\cdot dr\leq C_{11}||f||_{A^{p,q}_{\beta}}.
$$

Let $\displaystyle \frac{q}{p}>1, \chi_\gamma(z)=\frac{1}{(1-|z|)^{\frac{\gamma}{pq} }}, z\in U=\{|z|<1\}, T=\{|z|=1\}, 0<\gamma<\gamma_0.$

We show our theorem in the unit disk since it is typical and repetition of our arguments leads to the proof in bounded pseudoconvex domains. We have repeating arguments above

$$
\int\limits_0^1(1-\rho)^\beta\left(\int\limits_T|T_\alpha(f)(\rho z)|^p\cdot d\xi dr\right)^{q/p}\rho d\rho\leq 
$$

$$
\leq C_{12}\int\limits_0^1(1-\rho)^\beta\left(\int\limits_0^1\frac{(1-r)^{\alpha p+2(p-1)}}{(1-r\rho)^{(\alpha+2)p-1}}             
\int\limits_{T}|f(r\xi)|^pd\xi dr\right)^{q/p}\cdot \rho d\rho=J(f).
$$ 

Using Holder's inequality and changing the  order of integral we have that 

$$
J(f) \leq C_{13}\int\limits_0^1\left(\int\limits_0^1\frac{(1-r)^{\alpha p+2(p-1)}}{(1-r\rho)^{(\alpha+2)p-1}\cdot (\chi_\gamma^{q/p}(r))} \right)\cdot      
\left(\int\limits_{T}|f(r\xi)|^pd\xi \right)^{q/p}\cdot
$$

$$
\cdot\left(\int\limits_0^1\frac{(1-r)^{\alpha p+2(p-1)}\cdot \chi_\gamma^{\frac{q}{q-p}}(r)}{(1-r\rho)^{(\alpha+2)p-1}}    dr\right) ^{\left(\frac{q}{q-p}\right)^{-1}\cdot (q/p)}(1-\rho)^{\beta} \leq
$$ 

$$
\leq C_{14}\int\limits_0^1 (1-\rho)^\beta\cdot \chi_\gamma^{q/p}(\rho)\cdot\left(\int\limits_0^1\frac{(1-r)^{\alpha p+2(p-1)}rdr}{(1-r\rho)^{(\alpha+2)p-1}\chi_\gamma^{q/p}(r)}    dr\right) \times\left(\int\limits_T|f(r\xi)|^pd\xi\right)^{q/p}dr\rho d\rho\leq
$$ 

$$
\leq C_{15}\int\limits_0^1\frac{(1-r)^{\alpha p+2(p-1)}}{\chi_\gamma^{q/p}(r)}  \cdot \left(\int\limits_T|f(r\xi)|^pd\xi\right)^{q/p} \left(\int\limits_0^1\frac{(1-\rho)^\beta\chi_\gamma^{q/p}(\rho) d\rho}{(1-r\rho)^{(\alpha+2)p-1}}\right) dr \leq
$$ 

$$
\leq C_{16}\int\limits_0^1\frac{(1-r)^{\alpha p+2(p-1)}\times(1-r)^\beta   (\chi_\gamma^{q/p}(r))}{(\chi_\gamma^{q/p}(r))(1-r)^{(\alpha+2)p-2}}  \cdot 
\left(\int\limits_T|f(r\xi)|^pd\xi\right)^{q/p} dr \leq C_2||f||_{A^{p,q}_{\beta}}.
$$ 

Theorem is proved.

{\bf Remark 3.}   We suppose that our theorem is valid for $A^{p,q}$ spaces when both $p, q\geq 1$ and even when $p$  is infinite.

We define mixed norm analytic function spaces in bounded strongly pseudoconvex  domains with  boundary (in product domains).

Let 
$$
A^{p_1,p_2, \dots ,p_n}_{\alpha_1, \alpha_2, \dots ,\alpha_n}(D^m)=
$$
$$
=\left\{f\in H(D\times \ldots \times D): \left(\int\limits_D \ldots
\left(\int\limits_D|f(z_1,\ldots, z_m)|^{p_1}dv_{\alpha_1}(z_1)
\right)^{p_2/p_1}\cdots dv_{\alpha_m}(z_m)\right)^{1/p_m}<\infty\right\}.
$$

Changing $A$ by $L$ we denote the larger space consisting of all measurable functions in  $D\times \ldots \times D$, where  $0<p_i<\infty, -1<\alpha_i<\infty,i =1,\ldots,m$.

{\bf Theorem 2. } 

Let $p_j>1, $ $i =1,\ldots,m,$ $\alpha_j>-1,$ $j=1,\ldots,m.$ Then we have that the $V_{\vec{\beta}}$ operator, where
$$
\left(V_{\vec{\beta}}\right)(f)=\int\limits_D\ldots\int\limits_D f(z_1,\ldots,z_m)\prod\limits^{m}_{j=1}K_{\beta_j+n+1}(z_j,w_j)dv_{\beta_1}(z_1)\ldots dv_{\beta_m}(z_m),
$$
for all $\beta_j>\beta_0,  j=1,\ldots,m$ $\beta_0,$ is large enough maps 
$L^{p_1,\ldots,p_n}_{\alpha_1, \ldots, \alpha_n}$ into $A^{p_1,\ldots,p_n}_{\alpha_1, \ldots, \alpha_n}$ space.

{\bf Remark 4.}  In the unit disk for $m=1$ this result is classical and well-known fact  (Bergman projection theorem). The proof uses only Forelly-Rudin estimate from Lemma 1 and Holders and Minkowski inequalities and $m=2$ and unit disk case is typical. We have in the unit disk $U=\{|z|<1\}, m=2$ 
case the following estimates. 

We denote as usual by $dm_2$ the normalized Lebegues measure in the unit disk $U$.

Put first 
$$
D_{\alpha_j}(\xi_j,z_j)=\frac{\alpha_j+1}{\pi}\cdot\frac{(1-|\xi_j|)^{\alpha_j}}{(1-\bar{\xi}_jz_j)^{\alpha_j+2}}, j=1,2;
$$
$$D_\alpha(z.\xi)=D_{\alpha_1}(z_1.\xi_1)\times D_{\alpha_2}(z_2.\xi_2);$$
$$
\frac{1}{p_j}+\frac{1}{q_j}=1, j=1,2, \chi_\gamma(z_1,z_2)=(1-|z_1|)^{\frac{-\gamma}{p_1q_1}}\cdot(1-|z_2|)^{\frac{-\gamma}{p_2q_2}}, z_j\in U, \xi_j\in U, j=1,2.
$$
$$
||V_{\vec{\beta}}f(\cdot,z_2)||_{A^{p}_{\alpha_1}}=\left(\int\limits_{U}|F(z_1,z_2)|^{p_1}(1-|z_1|)^{\alpha_1}dm_2(z_1)\right)^{1/p_1}\leq
$$
$$
\leq C_{17}\left(\int\limits_{U}\left(\int\limits_{U}|D_{\alpha_1}(z_1,\xi_1)|\cdot|D_{\alpha_2}(z_2,\xi_2)|\cdot|f(\xi_1,\xi_2)|dm_2(\xi_1)dm_2(\xi_2)\right)^{p_1} (1-|z_1|)^{\alpha_1}dm_2(z_1)\right)\leq
$$
$$
\leq C_{18}\int\limits_{U}|D_{\alpha_2}(z_2,\xi_2)|\left(\int\limits_{U}\left(\int\limits_{U}|D_{\alpha_1}(z_1,\xi_1)|\cdot|f(\xi_1,\xi_2)|dm_2(\xi_1)\right)^{p_1} (1-|z_1|)^{\alpha_1}dm_2(z_1)\right)^{1/p_1} dm_2(\xi_2)\leq
$$
$$
\leq C_{19}\int\limits_{U}|D_{\alpha_2}(z_2,\xi_2)|\Bigg(\int\limits_{U}\left(\int\limits_{U}\frac{|D_{\alpha_1}(z_1,\xi_1)|\cdot|f(\xi_1,\xi_2)|^{p_1}dm_1(\xi_1)}{|\chi_\gamma^{p_1}(\xi_1,\xi_2)|}\right)\times
$$
$$\times\left(\int\limits_{U}|D_{\alpha_1}(z_1,\xi_1)|\cdot\chi_\gamma^{p_1}(\xi_1,\xi_2)dm_2(\xi_1)\right)^{p_1/q_1} (1-|z_1|)^{\alpha_1}dm_2(z_1)\Bigg)^{1/p_1} dm_2(\xi_2)\leq
$$
$$
\leq C_{20}\int\limits_{U}|D_{\alpha_2}(z_2,\xi_2)|\Bigg(\int\limits_{U}\frac{|f(\xi_1,\xi_2)|^{p}}{\chi_\gamma^{p_1}(\xi_1,\xi_2)} \int\limits_{U}|D_{\alpha_1}(z_1,\xi_1)|\chi_\gamma^{p_1}(z_1,\xi_2)(1-|z_1|)^{\alpha_1}dm_2(z_1)dm_2(\xi_2)\Bigg)^{1/p_1} dm_2(\xi_2)\leq
$$
$$
\leq C_{21}\int\limits_{U}|D_{\alpha_2}(z_2,\xi_2)|
\cdot ||f(\cdot, \xi_2)||
_{L^{p_1}_{\alpha_1}}dm_2(\xi_2).
$$

Then we have that
$$
||V_{\vec{\beta}}f(\cdot,z_2)||^{p_2}_{A^{p_1,p_2}_{\alpha_1,\alpha_2}}\leq C_{22}\int\limits_{U}\left(\int\limits_{U}|D_{\alpha_2}(z_1,\xi_2)|\cdot ||f(\cdot, \xi_2)||
_{L^{p_1}_{\alpha_1}}dm_2(\xi_2)\right)^{p_2}(1-|z_2|)^{\alpha_2}dm_2(z_2).
$$

Using again Holders inequality with $p_2$ and Forelly-Rudin estimate and changing the order of integration we get what we need

$$
||V_{\vec{\beta}}(f)||^{p_2}_{A^{p_1,p_2}_{\alpha_1,\alpha_2}}\leq C ||f||^{p_2}_{L^{p_1,p_2}_{\alpha_1,\alpha_2}}.
$$

Theorem is proved.

{\bf Remark 5.} We conjecture that similar results are valid when
one index is infinite in our mixed norm $A^{p_1,...,p_m}$ spaces.

Our projection theorems can be used in the study of action of Hankel and Toeplits operators in bounded strongly pseudoconvex domains (see for example [6]).

{\bf Information about the authors}\\

{\bf Romi F. Shamoyan}, PhD, Bryansk State University, Bryansk, Russia; \\rsham@mail.com.\\

{\bf Elena B. Tomashevskaya}, PhD, Bryansk State Tecnical University, Bryansk, Russia; tomele@mail.ru.

 \end{document}